\documentclass[12pt]{article}
\usepackage{amssymb}

\newcommand{\Z}{{\mathbb Z}}
\newcommand{\F}{{\mathbb F}}
\newcommand{\C}{{\mathbb C}}
\newcommand{\R}{{\mathbb R}}
\newcommand{\FF}{{\mathbb F}}
\newcommand{\eps}{\varepsilon}

\newtheorem{theorem}{Theorem}
\newtheorem{lemma}{Lemma}

\title{On sums and products in $\C[x]$}

\author{Ernie Croot and Derrick Hart}

\begin{document}

\maketitle

\section{Introduction}

Suppose that $S$ is a subset of a ring $R$ (in our case, the
real or complex numbers), and define
$$
S.S\ :=\ \{ st\ :\ s,t \in S\},\ {\rm and\ } S+S\ :=\ \{s+t\ :\ s,t \in S\}.
$$ 
An old problem of Erdos and Szemer\'edi \cite{erdos_sz} is to show
that 
$$
|S.S| + |S+S|\ \gg\ |S|^{2-o(1)}.
$$
Partial progress on this problem has been achieved by Erd\H os and
Szemer\'edi \cite{erdos_sz}, Nathanson \cite{nathanson}, 
Ford \cite{ford}, Elekes \cite{elekes}, 
and Solymosi in \cite{solymosi1} and in the astounding article 
\cite{solymosi2}.  There is the work of Bourgain-Katz-Tao \cite{BKT} 
and Bouragin-Glibichuk-Konyagin \cite{BGK} extending the results to
$\FF_p$, and then also the recent work of Tao \cite{tao} which extends
these results to arbitrary rings.
\bigskip

We will prove a theorem below (Theorem \ref{main_theorem2}), 
which holds for the polynomial ring $\C[x]$ {\it unconditionally}; 
but, under the assumption of a certain $24$-term version of Fermat's Last 
Theorem, it holds for $\Z_+$, the
positive integers, as we will explain below.  This {\it unconditional} 
result for $\C[x]$, for which there presently is no 
analogue for $\C$, $\R$, $\Z$ or $\F_p$
(though, M.-C. Chang has a related sort of result for 
$\Z$ as we discuss in remarks below.), is as follows.

\begin{theorem} \label{main_theorem2} There exists an absolute constant
$c > 0$ such that the following holds for all sufficiently large
sets $S$ of monic polynomials of $\C[x]$ of size $n$:
$$
|S.S|\ <\ n^{1+c}\ \ \Longrightarrow\ \ |S+S|\ \gg\ n^2.
$$
\end{theorem}

\noindent {\bf Remark 1.}  M.-C. Chang \cite{mc1} 
has shown that if $S$ is a set of integers such that 
$|S.S| < n^{1+\eps}$, then $|S+S| > n^{2-f(\eps)}$, where $f(\eps) \to 0$
as $\eps \to 0$.
\bigskip

\noindent {\bf Remark 2.}  It is perhaps possible to replace the 
condition that the polynomials be ``monic'' with the condition that none
is a scalar multiple of another; however, it will make the proof more
complicated (if our method of proof is used).
\bigskip

Now let us consider the following conjecture, which can be thought of as
a $24$-term extension of Fermat's Last Theorem.
\bigskip

\noindent {\bf Conjecture.}  For all $m \geq 1$ sufficiently large, the only
solutions to the Diophantine equation
\begin{equation} \label{epsm_equation}
\eps_1 x_1^m + \eps_2 x_2^m + \cdots + \eps_{24} x_{24}^m\ =\ 0,\ 
\eps_i = \pm 1,\ x_1,...,x_{24} \in \Z_+
\end{equation}
are the trivial ones; that is, solutions where each $\eps_i x_i^m$ 
can be paired with its negative $\eps_j x_j^m$, so that
$x_i = x_j$ and $\eps_i = -\eps_j$.  
\bigskip

Under the assumption of this conjecture, we have that 
Theorem \ref{main_theorem2} holds for when $S$ is a set of positive
integers (instead of monic polynomials); in other words, 

\begin{theorem} \label{main_theorem3}  Suppose that the above {\bf Conjecture}
is true.  Then, there exists an absolute constant
$c > 0$ such that the following holds for all sufficiently large
sets $S$ of $n$ positive integers:
$$
|S.S|\ <\ n^{1+c}\ \ \Longrightarrow\ \ |S+S|\ \gg\ n^2.
$$
\end{theorem}

Because the proof of this theorem is virtually identical to the proof of
Theorem \ref{main_theorem2}, where every use of Theorem \ref{poly_theorem}
below is simply replaced with the above {\bf Conjecture}, we simply omit
the proof. 
\bigskip

We also prove the following theorem, which extends a result of  
Bourgain and Chang \cite{bourgain} from $\Z$ to the ring $\C[x]$.

\begin{theorem} \label{main_theorem}  Given a real $c \geq 1$ and
integers $\ell, k \geq 1$, the following holds for all $n$ sufficiently large:
Suppose that $S$ is a set of $n$ polynomials of $\C[x]$, 
where none is a scalar multiple of another, and suppose that 
$$
|S^\ell|\ =\ |S.S...S|\ <\ n^c.
$$
Then,
$$
|kS|\ >\ n^{k-f(c,k,\ell)},
$$
where for fixed $c$ and $k$, $f(c,k,\ell) \to 0$ as $\ell \to \infty$.
\end{theorem}

\noindent {\bf Remark 1.}  It is possible to generalize our method of proof to 
the case of polynomials over $\F_p[x]$, but in 
that context there are thorny issues concerning the vanishing of 
certain Wronskian determinants that make the problem difficult.  There
are also issues that come up in handling $p$th powers of polynomials,
which the differentiation mapping sends to the $0$ polynomial.
\bigskip

\noindent {\bf Remark 2.}  If the analogue of this theorem for when
$S \subseteq \C$ could be proved, then it would also provide a proof 
to the above, since we could locate an element $\alpha \in \C$ such that
the evaluation map $\alpha : f \in \C[x] \to f(\alpha) \in \C$
preserves the structure of the sums and products of these polynomials.  
\bigskip

\noindent {\bf Remark 3.}  Using the above results, one can easily prove
analogues of them for the polynomial ring $\C[x_1,x_2,...,x_k]$, simply
by applying an evaluation map $\psi : \C[x_1,...,x_k] \to 
\C[x_1,\alpha_1,...,\alpha_k]$, for some carefully chosen 
$\alpha_1,...,\alpha_k$.
\bigskip

The proofs of both of the above theorems rely on the following basic
fact about polynomials, using ideas of 
Mason \cite{mason}, which we prove near the end of the paper.

\begin{theorem} \label{poly_theorem}  
For every $k \geq 2$, there exists an exponent $M$ such that
there are no polynomial solutions to 
\begin{equation} \label{poly_equation}
f_1(x)^m + \cdots + f_k(x)^m\ =\ 0,\ m \geq M,
\end{equation}
where no polynomial is a scalar multiple of another.
\end{theorem}

\noindent {\bf Remark.}  Previously, it had been proved that such 
an equation has no solutions under
the condition that all the polynomials are pairwise coprime.
\bigskip

We will also make use of the Ruzsa-Plunnecke \cite{plunnecke},
\cite{ruzsa} inequality, stated as follows.

\begin{theorem} \label{ruzsa_plunnecke} 
Suppose that $S$ is a finite subset of an additive abelian group, and that
$$
|S+S|\ \leq\ K|S|.
$$
Then,
$$
|kS - \ell S|\ =\ |S+S\cdots + S - S - S - \cdots - S|\ \leq\ K^{k+\ell} |S|.
$$
\end{theorem}

\section{Proof of Theorem \ref{main_theorem2}}

By invoking Theorem \ref{poly_theorem}, let $M \geq 1$ be the smallest
value such that the polynomial equation
$$
f_1(x)^m + \cdots + f_k(x)^m\ =\ 0
$$
has no solutions for $m \geq M$, and $k \leq 24$, assuming no $f_i$
is a constant multiple of another.  And, when we generalize the
present proof (of Theorem \ref{main_theorem2}) to handle the proof of 
Theorem \ref{main_theorem3}, we just assume that $M$ is such that if 
$m \geq M$, then (\ref{epsm_equation}) has no non-trivial solutions. 

We assume throughout that 
$$
|S.S|\ <\ n^{1+c},
$$
where $c > 0$ is some parameter that can be determined by working through
the proof -- the point is that, although $c$ will be quite small,
it will be possible to take it to be some explicit value.

Let 
$$
\eps\ =\ c(M+1),
$$
and note that by the Ruzsa-Plunnecke inequality 
(Theorem \ref{ruzsa_plunnecke}) we have that since $|S.S| < n^{1+c}$,
$$
|S^{M+1}|\ =\ |S...S|\ <\ n^{1+\eps}.
$$
(The sole use of $\eps > 0$ in the rest of the proof is to simplify
certain expressions.)
\bigskip

Let us begin by supposing that 
$$
|S+S|\ =\ o(n^2).
$$
Then, for all but at most $o(n^2)$ pairs $(x_1,x_2) \in S \times S$,
there exists $(x_3,x_4)$ with
$$
\{x_3, x_4\}\ \neq\ \{x_1, x_2\},
$$
such that
$$
x_1 + x_2\ =\ x_3 + x_4.
$$
Let $P$ denote the set of all such $n^2 - o(n^2)$ pairs $(x_1,x_2)$.
It is clear that there exists a bijection 
$$
\varphi\ :\ P\ \to\ P,
$$
where $\varphi$ maps pairs having sum $s$ to pairs having sum $s$,
and yet where if 
$$
(x_3, x_4)\ =\ \varphi((x_1,x_2)),
$$
then
$$
\{x_3, x_4\}\ \neq\ \{x_1,x_2\}.
$$

Using such a pairing $\varphi$ of pairs $(x_1,x_2)$, we then define
a set of quadruples
$$
Q\ :=\ 
\{ (x_1,x_2,x_3,x_4)\ :\ (x_1,x_2) \in P,\ (x_3,x_4) = \varphi((x_1,x_2))\}.
$$
Note that
$$
|Q|\ =\ |P|\ \sim\ n^2,
$$
and
$$
(x_1,x_2,x_3,x_4) \in Q\ \ \Longrightarrow\ \ x_1+x_2-x_3-x_4\ =\ 0.
$$

\subsection{A lemma about quadruples}

To proceed further we require the following lemma.

\begin{lemma}  There exists a 5-tuple $a,b,c,d,t \in S$ such that for 
$\gg n^{2-4\eps}$ quadruples 
$$
(x_1,x_2,x_3,x_4)\ \in\ Q
$$
we have that 
$$
t^M(x_1,x_2,x_3,x_4)\ =\ (a t_1^M, bt_2^M, c t_3^M, d t_4^M),\ 
{\rm where\ } t_1,t_2,t_3,t_4 \in S.
$$
\end{lemma} 

\noindent {\bf Proof of the lemma.}  Let $N$ denote the number of pairs 
$(x_1,t) \in S \times S$ for which there are fewer than 
$n^{1-\eps}/40$ pairs $(a,t_1) \in S \times S$ satisfying 
$$
x_1 t^M\ =\ a t_1^M.
$$
Clearly,
\begin{eqnarray*}
N\ &\leq&\ (n^{1-\eps}/40) |\{ \alpha \beta^M\ :\ \alpha, \beta \in S\}| \\
&\leq&\ (n^{1-\eps}/40) |S^{M+1}| \\
&<&\ n^2/40.
\end{eqnarray*} 

Given $x_1 \in S$ we say that $t$ is ``bad'' if $(x_1,t)$ is one of the 
pairs counted by $N$; otherwise, we say that $t$ is ``good''.  
From the bounds above, it is clear that 
more than $4n/5$ values $x_1 \in S$ have $\leq n/8$ bad values of 
$t$; for, if there were fewer than $4n/5$ such $x_1 \in S$, then 
$\geq n/5$ have $> n/8$ bad values of $t$, which would show that 
$N > n^2/40$, a contradiction.

It follows that more than $3n^2/5$ pairs $(x_1,x_2) \in S^2$ have 
the property that both $x_1$ and $x_2$ have at most $n/8$ bad values
of $t$.  So, there are at least 
$$
n\ -\ 2 n/8\ =\ 3n/4
$$
values of $t$ that are ``good'' for both $x_1$ and $x_2$. 

Clearly, then, by the pigeonhole principle, there are 
$$
\gtrsim\ n^2/5\ \ {\rm quadruples\ } (x_1,x_2,x_3,x_4) \in Q,
$$
such that there are at least $n/2$ values of $t$ that are 
good for all $x_1,x_2,x_3$ and $x_4$ at the same time (disguised in 
what we are doing here is the fact that 
$\varphi$ is a bijection from $P \to P$).  

When $t$ is good for all $x_1,x_2,x_3$ and $x_4$, we say that 
it is ``good'' for the quadruple $(x_1,x_2,x_3,x_4)$.

By the pigeonhole principle again, there exits $t \in S$ that is
good for at least 
$$
\gtrsim n^2/10
$$
quadruples of $Q$.  

Now suppose that $(x_1,x_2,x_3,x_4) \in Q$ is one of these $\sim n^2/10$
quadruples, for some particular value of $t$.  We define
\begin{eqnarray*} 
A\ &:=&\ \{ a \in S\ :\ a t_1^M = x_1 t^M,\ t_1 \in S\} \\
B\ &:=&\ \{ b \in S\ :\ b t_2^M = x_2 t^M,\ t_2 \in S\} \\
C\ &:=&\ \{ c \in S\ :\ c t_3^M = x_3 t^M,\ t_3 \in S\} \\
D\ &:=&\ \{ d \in S\ :\ d t_4^M = x_4 t^M,\ t_4 \in S\}.
\end{eqnarray*}

Under the assumption that $t$ is good, we have 
$$
|A|, |B|, |C|, |D|\ \geq\ n^{1-\eps}/40.
$$
So, among the $\gtrsim n^2/10$ quadruples for which $t$ is good, by the
pigeonhole principle again, there exists $a,b,c,d\ \in\ S$ for which
there are 
$$
\gtrsim\ (n^2/10) (n^{-\eps}/40)^4\ \gg\ n^{2-4\eps}
$$
quadruples $(x_1,x_2,x_3,x_4) \in Q$ satisfying
$$
x_1 t^M = a t_1^M,\ ...,\ x_4 t^M = d t_4^M.
$$
This completes the proof of the lemma.
\hfill $\blacksquare$

\subsection{Resumption of the proof of Theorem \ref{main_theorem2}}

Upon applying the previous lemma, let $Q'$ denote the set of 
quadruples
$$
(t_1,t_2,t_3,t_4)\ \in\ S^4,
$$
such that
$$
x_1 t^M = a t_1^M,\ ...,\ x_4 t^M = d t_4^M,
$$
where $(x_1,x_2,x_3,x_4)$ is one of the $\gg n^{2-4\eps}$ quadruples
satisfying the conclusion of the lemma.  Note that
$$
a t_1^M + b t_2^M - c t_3^M - d t_4^M\ =\ 0.
$$
\bigskip

Now we claim that we can find three quadruples
$$
(t_1,...,t_4),\ (u_1,...,u_4),\ (v_1,...,v_4)\ \in\ Q',
$$
such that
\begin{equation} \label{no1}
t_2/t_1 \neq u_2/u_1,\ t_2/t_1 \neq v_2/v_1,\ u_2/u_1 \neq v_2/v_1,
\end{equation}
and
\begin{equation} \label{no2}
t_4/t_3 \neq u_4/u_3,\ t_4/t_3 \neq v_4/v_3,\ u_4/u_3 \neq v_4/v_3.
\end{equation}

The reason that we can find such quadruples is as follows:  Suppose
we fix $(t_1,...,t_4)$ to be any quadruple of $Q'$, and suppose that
we pick $(u_1,...,u_4)$ in order to attempt to avoid
\begin{equation} \label{yes1}
t_2/t_1 = u_2/u_1.
\end{equation}
Let $r = t_2/t_1 \in S/S$.  
If (\ref{yes1}) holds for all $(u_1,u_2)$, it means that
each pair $(u_1,u_2)$ must be of the form
$$
(u_1,u_2) = r (t_1,t_2).
$$
But, by Ruzsa-Plunnecke, 
$$
|S/S|\ <\ n^{1+\eps},
$$
so there are at most $n^{1+\eps}$ pairs among the first two 
entries of quadruples of $Q'$; but, since there at $\gg n^{2-4\eps}$ 
quadruples of $Q'$, and since the first two coordinates of the quadruple
determine the second pair of coordinates (via the mapping $\varphi$) 
we clearly must not have that all pairs $(u_1,u_2)$ satisfy
(\ref{yes1}).  In fact, there are $\gg n^{2-\eps} - O(n^{1+\eps})$ 
pairs (quadruples) to choose from!

In a similar vein, we can pick 
$(u_1,...,u_4)$ and $(v_1,...,v_4)$, so that all the remaining conditions 
(\ref{no1}) and (\ref{no2}) hold.
\bigskip

Before proceeding with the rest of the proof, it is worth pointing out 
that the sort of condition on our four-tuples that we cannot so easily
force to hold is, for example, 
$$
t_1/t_3 \neq u_1/u_3.
$$
The reason for this is that we do not have a good handle on how many
pairs $(t_1,t_3)$ or $(u_1,u_3)$ there are -- $t_3$ (or $u_3$) may
be related to $t_1$ (or $u_1$) in a completely trivial way, leading
to few pairs.

\subsection{A lemma about submatrices}

Now we require a lemma concerning the matrix
$$
T\ :=\ \left [ \begin{array}{cccc} t_1^M & t_2^M & t_3^M & t_4^M \\ 
u_1^M & u_2^M & u_3^M & u_4^M \\
v_1^M & v_2^M & v_3^M & v_4^M \end{array} \right ].
$$

\begin{lemma} \label{nonsingular_lemma} 
Every $3 \times 3$ submatrix of $T$ is non-singular.
\end{lemma}

\noindent {\bf Proof of the lemma.}  We show that $3 \times 3$
matrices are non-singular via contradiction:  Suppose, on the contrary, 
that some $3 \times 3$ submatrix {\it is} singular, and without
loss assume (that it is the `first' $3 \times 3$ submatrix): 
$$
\left | \begin{array}{ccc} t_1^M & t_2^M & t_3^M \\ u_1^M & u_2^M & u_3^M \\
v_1^M & v_2^M & v_3^M \end{array} \right |\ =\ 0.
$$
Expanding out the determinant into a polynomial in its entries, 
we see that it produces a sum of $M$th powers equal to $0$.
Furthermore, since all the $t_i,u_j,v_k$ are monic, none is a scalar multiple
of another, except for factors $\pm 1$.  Since no non-trivial sum
of 24 or fewer $M$th powers of polynomials can equal $0$, it follows that
each must be matched with its negative, in order for this sum of $M$th
powers to equal $0$.

Note, then, that there are $6 = 3!$ possible matchings that can produce a 
$0$ sum.  Consider now the matching (assuming we have taken $M$th roots)
\begin{eqnarray*}
t_1 u_2 v_3\ &=&\ t_3 u_2 v_1 \\
t_2 u_3 v_1\ &=&\ t_2 u_1 v_3 \\
t_3 u_1 v_2\ &=&\ t_1 u_3 v_2.
\end{eqnarray*}
This matching implies
\begin{equation} \label{useme}
t_3/t_1\ =\ u_3/u_1\ =\ v_3/v_1.
\end{equation}
Some of the other possible matchings will lead to equations such
as $t_1/t_2 = u_1/u_2$, which we have said was impossible by design;
but, there is one other viable chain of equations that we get, using
one of these matchings, and that is 
$$
t_3/t_2\ =\ u_3/u_2\ =\ v_3/v_2.
$$
In what follows, whether this chain holds, or (\ref{useme}) holds 
makes little difference, so we will assume without loss of generality
that (\ref{useme}) holds, and will let $r = t_3/t_1$ denote the
common ratio, which we note is a rational function.

It is clear that assuming that this matching holds, we can reduce  
the equations
\begin{eqnarray*}
at_1^M + b t_2^M - c t_3^M - d t_4^M\ &=&\ 0 \\
au_1^M + b u_2^M - c u_3^M - d u_4^M\ &=&\ 0 \\
av_1^M + b v_2^M - c v_3^M - d v_4^M\ &=&\ 0
\end{eqnarray*}
to
$$
\left [ \begin{array}{ccc} t_1^M & t_2^M & t_4^M \\
u_1^M & u_2^M & u_4^M \\
v_1^M & v_2^M & v_4^M \end{array} \right ] 
\left [ \begin{array}{c} a - cr^M \\ b \\ -d \end{array} \right ]
\ =\ \left [ \begin{array}{c} 0 \\ 0 \\ 0 \end{array} \right ].
$$
Since all the elements of $S$ are monic, none can be $0$, and so this
column vector is not the $0$ vector.  It follows that the $3 \times 3$
matrix here is singular (from the fact that one of the $3 \times 3$ submatrices
of $T$ is singular, we just got that another was singular).  Upon expanding
the determinant of this matrix into a polynomial of its entries, 
in order to get it to be $0$ we must have a matching, much like the
one that produced (\ref{useme}).  As before, we will get two viable chains
of equations:  Either
\begin{equation} \label{first_chain}
t_4/t_1\ =\ u_4/u_1\ =\ v_4/v_1,
\end{equation}
or
\begin{equation} \label{second_chain}
t_4/t_2\ =\ u_4/u_2\ =\ v_4/v_2.
\end{equation}
Let us suppose that (\ref{first_chain}) holds, and let $s$ be the common ratio,
which is of course a rational function.
Then, it follows that
$$
\left [ \begin{array}{cc} t_1^M & t_2^M \\ u_1^M & u_2^M \end{array} 
\right ] \left [ \begin{array}{c} a - cr^M - ds^M \\ b \end{array} \right ]
\ =\ \left [ \begin{array}{c} 0 \\ 0 \end{array} \right ].
$$
Since $(t_1^M,t_2^M)$ and $(u_1^M,u_2^M)$ are independent, this matrix
is non-singular; so, the column vector here must be the $0$ vector,
which is impossible since $b \neq 0$.  

So, (\ref{second_chain}) above must hold.  Redefining $s$ to be the common
ratio $t_4/t_2$ here, it leads to the equation
$$
\left [ \begin{array}{cc} t_1^M & t_2^M \\ u_1^M & u_2^M \end{array}
\right ] \left [ \begin{array}{c} a - cr^M \\ b - ds^M \end{array} \right ]
\ =\ \left [ \begin{array}{c} 0 \\ 0 \end{array} \right ].
$$
Since the matrix is non-singular, it follows that
$$
a = cr^M\ \ {\rm and\ \ } b = ds^M.
$$
But this implies that
\begin{equation} \label{bad1}
a t_1^M\ =\ c r^M t_1^M\ =\ c t_3^M,
\end{equation}
and
\begin{equation} \label{bad2}
b t_2^M\ =\ d s^M t_2^M\ =\ d t_4^M,
\end{equation}
So, $(t_1,t_2,t_3,t_4)$ couldn't have been a quadruple of $Q'$, because
if so, then for some $(x_1,x_2,x_3,x_4) \in Q$ we would have had
$$
t^M(x_1,x_2,x_3,x_4)\ =\ (a t_1^M, b t_2^M, c t_3^M, d t_4^M),
$$
and then the equations (\ref{bad1}) and (\ref{bad2}) imply that
$$
(x_1,x_2) = (x_3,x_4),
$$
a contradiction.  This completes the proof of the lemma.
\hfill $\blacksquare$

\subsection{Continuation of the proof}

In addition to the quadruples $(t_1,...,t_4), (u_1,...,u_4), (v_1,...,v_4)$,
let $(w_1,...,w_4)$ be another quadruple (we will later show that
we can choose it in such a way that we contradict the assumption 
$|S+S| = o(n^2)$).  We have that
$$
\Gamma\ :=\ \left [ \begin{array}{cccc} t_1^M & t_2^M & t_3^M & t_4^M \\
u_1^M & u_2^M & u_3^M & u_4^M \\
v_1^M & v_2^M & v_3^M & v_4^M \\
w_1^M & w_2^M & w_3^M & w_4^M \end{array} \right ] 
$$
is singular, as the vector $(a,b,-c,-d)$ (written as a column vector)
is in its kernel.  
Expanding out its determinant, we find that it must be $0$; and,
we know from Theorem \ref{poly_theorem} that this is impossible, except
if we can match up terms, as we did in the proof of Lemma 
\ref{nonsingular_lemma}.  Just so the reader is clear, an example of
just one equation from such a matching is perhaps, say, 
$$
(t_1 u_2 v_3 w_4)^M\ =\ (t_2 u_1 v_3 w_4)^M.
$$ 

In total, there will be $12 = 4!/2$ different equations that make up such 
a matching.  Let us suppose that in a hypothetical matching we got
an equation of the form
$$
(t_i u_j v_k w_1)^M\ =\ (t_{i'} u_{j'} v_{k'} w_2)^M,
$$
and for our purposes we just need to write this as
\begin{equation} \label{badalpha}
\alpha w_1^M\ =\ \beta w_2^M,
\end{equation}
where $(\alpha, \beta) \in \C[x] \times \C[x]$ can be any of at most 
$6^2$ polynomials (according to which combination of $i,j,k,i',j',k'$
is chosen).  So, we would have $w_2^M/w_1^M = \alpha/\beta$; that is,
$w_2/w_1$ takes on at most $36$ possible values.  What this would mean
is that the pair $(w_1,w_2)$ is essentially determined by $w_1$,
and that it can take on at most $36 n$ possible values -- far too
few to consume the majority of the $\gg n^{2-4\eps}$ pairs $(w_1,w_2)$
making up the first part of a quadruple of $Q'$.  

We conclude that all but $O(n)$ of the quadruples $(w_1,w_2,w_3,w_4) \in Q'$
cannot lead to a solution to (\ref{badalpha}) under any matching; 
moreover, all but $O(n)$ will also avoid
\begin{equation} \label{situation}
\alpha w_3^M\ =\ \beta w_4^M.
\end{equation}

We also have that there are at most $O(n)$ quadruples can lead to solutions 
to 
\begin{equation} \label{nopair1}
\alpha w_1^M\ =\ \beta w_3^M,\ {\rm and\ } \alpha' w_1^M\ =\ \beta' w_4^M,
\end{equation}
because it would imply that $w_4/w_3$ is fixed, and we are back in
the situation (\ref{situation}).  Furthermore, there are at most 
$O(n)$ quadruples leading to solutions to any of the following pairs:
\begin{equation} \label{nopair2}
\alpha w_2^M\ =\ \beta w_3^M,\ {\rm and\ } \alpha' w_2^M\ =\ \beta' w_4^M,
\end{equation}
or
\begin{equation} \label{nopair3}
\alpha w_3^M\ =\ \beta w_1^M,\ {\rm and\ } \alpha w_3^M\ =\ \beta w_2^M,
\end{equation}
or
\begin{equation} \label{nopair4}
\alpha w_4^M\ =\ \beta w_1^M,\ {\rm and\ } \alpha w_4^M\ =\ \beta w_2^M.
\end{equation}
\bigskip

We can also avoid a matching that produces three equations (indexed by
$j$) of the form
\begin{equation} \label{atriple}
\alpha_j w_i^M = \beta_j w_i^M,\ j=1,2,3,
\end{equation}
because it would mean that the $3 \times 3$ submatrix of $T$ with the
$i$th column deleted, is singular (these $\alpha_j,\beta_j$ have the
property that the determinant of this submatrix is
$\sum_{j=1}^3 (\alpha_j - \beta_j)$).

Furthermore, we cannot even have a {\it pair} of equations of the type 
(\ref{atriple}), for the same value $i$, 
because it would imply that in fact we get three equations
upon taking a product of the two and doing some cancellation; for example,
suppose that $i=4$, and that we have two equations of (\ref{atriple})
holding.  Then, there is a matching between two pairs of terms,
upon expanding the determinant of the following matrix in terms of its 
entries: 
\begin{equation} \label{matrix_focus}
\left [ \begin{array}{ccc} t_1^M & t_2^M & t_3^M \\
u_1^M & u_2^M & u_3^M \\
v_1^M & v_2^M & v_3^M \end{array} \right ].
\end{equation}
The matching corresponds, say, to $\alpha_1 = \beta_1$ and
$\alpha_2 = \beta_2$; and, say, these correspond to the equations 
\begin{eqnarray*}
(t_1 u_2 v_3)^M\ &=&\ (t_3 u_2 v_1)^M \\ 
(t_2 u_3 v_1)^M\ &=&\ (t_2 u_1 v_3)^M.
\end{eqnarray*}
Multiplying left and right sides together, and cancelling, 
produces
$$
(t_1 u_3)^M\ =\  (t_3 u_1)^M;
$$ 
and, multiplying both sides by $v_2^M$ produces the missing matching
$$
(t_1 u_3 v_2)^M\ =\ (t_3 u_1 v_2)^M,
$$
which proves that the matrix (\ref{matrix_focus}) is singular.  
But this is impossible, since it contradicts Lemma \ref{nonsingular_lemma}.
We conclude therefore that, as claimed, we cannot have even a pair
of equations from (\ref{atriple}) hold, for any $i=1,...,4$.
\bigskip

We have eliminated a great many possible matchings that could occur in
order that the matrix $\Gamma$ be singular:  Our $12$ equations in 
a matching can include at most one of each the four types
$$
\alpha w_1^M = \beta w_1^M,\ \alpha' w_2^M = \beta' w_2^M,\ 
\alpha'' w_3^M = \beta'' w_3^M,\ \alpha''' w_4^M = \beta''' w_4^M,
$$
and so must include at least $8$ of the form
\begin{equation} \label{theeight}
\alpha w_1^M = \beta w_3^M,\ \alpha' w_1^M = \beta' w_4^M, 
\alpha'' w_2^M = \beta'' w_3^M,\ \alpha''' w_2^M = \beta''' w_4^M,
\end{equation}
all the while avoiding pairs (\ref{nopair1}), (\ref{nopair2}),
(\ref{nopair3}), and (\ref{nopair4}).
\bigskip

Let us consider what would happen if we had a pair
\begin{equation} \label{notmany}
\alpha w_1^M = \beta w_3^M\ \ {\rm and\ \ } \alpha''' w_2^M = \beta''' w_4^M,
\end{equation}
or a pair
\begin{equation} \label{notmany2}
\alpha' w_1^M = \beta' w_4^M\ \ {\rm and\ \ } \alpha'' w_2^M = \beta'' w_3^M,
\end{equation}
both holding.  Without loss of generality in what follows, we
just assume that (\ref{notmany}) holds.  Then, we would have that the equation
$$
a w_1^M + b w_2^M - c w_3^M - d w_4^M\ =\ 0
$$
becomes
$$
(a - c \alpha/\beta) w_1^M + (b - d \alpha'''/\beta''') w_2^M\ =\ 0.
$$ 
Now, if we had any other quadruple $(z_1,z_2,z_3,z_4) \in Q'$
that also satisfied 
$$
\alpha z_1^M = \beta z_3^M\ \ {\rm and\ \ } \alpha''' z_2^M = \beta''' z_4^M,
$$
we would likewise get
$$
(a - c \alpha/\beta) z_1^M + (b - d \alpha'''/\beta''') z_2^M\ =\ 0,
$$
and then we would get the equation
$$
\left [ \begin{array}{cc} w_1^M & w_2^M \\ z_1^M & z_2^M \end{array} \right ]
\left [ \begin{array}{c} a - c\alpha/\beta \\ b - d \alpha'''/\beta''' 
\end{array} \right ]\ =\ \left [ \begin{array}{c} 0 \\ 0 \end{array} \right ].
$$
So, either we get 
$$
a = c \alpha/\beta,\ b = d \alpha'''/\beta''',
$$
or else this matrix is singular.  If the former holds, it implies that
$$
(a w_1^M, bw_2^M)\ =\ (c w_3^M, d w_3^M),
$$
which contradicts the hypotheses about the sets of quadruples $Q$ and $Q'$.

So, the matrix must be singular; in other words,
$$
(w_1, w_2)\ =\ \gamma (z_1,z_2),\ \gamma \in S/S.
$$
But since by the Ruzsa-Plunnecke inequality $|S/S| < n^{1+\eps}$, 
there can be at most $n^{1+\eps}$ such vectors $(z_1,z_2)$, given
$(w_1,w_2)$.

In particular, this means that there can be only very few quadruples 
$(w_1,w_2,w_3,w_4) \in Q'$ that satisfy (\ref{notmany}) or 
(\ref{notmany2}), for 
any particular combination of $\alpha, \beta, \alpha'', \beta'', 
\alpha''', \beta'''$.
Since there are at most $6^6$ possibilities for all these, 
(they are products of entries from $T$), we deduce
that there are $\gg n^{2-4\eps} - O(n^{1+\eps})$ quadruples of $Q'$
that do not satisfy a pair of equations of the sort 
(\ref{notmany}) or (\ref{notmany2}).  
So, we may safely assume that $(w_1,...,w_4)$ does not satisfy
these equations.
\bigskip

We may assume, then, that all 8 of our equations of the type 
(\ref{theeight}) involve the same pair of $w_i$'s, and do not
involve $w_1,w_2$ or $w_3,w_4$.  So, for example, our matching includes
8 equations of the form 
$$
\alpha_j w_1^M = \beta_j w_3^M,\ j=1,...,8
$$
(or 8 analogous equations for $w_1,w_4$ or $w_2,w_3$ or $w_2,w_4$). 
But, thinking about where such equations come from (from a matching
on the matrix $\Gamma$), there simply {\it cannot be} 8 equations:
$\alpha_j$ is any of the 6 terms (times possibly $-1$) making up the
determinant of the submatrix of $\Gamma$ gotten by deleting the 
last row and the first column (or just the first column of $T$), 
while the $\beta_j$ corresponds to possible
terms in the determinant of the submatrix gotten by deleting the last
row and third column of $\Gamma$.  So, each of the six terms in the first
determinant, matched with a unique term of the second, produces only
6 equations, not 8.  

We have now exhausted all of the possibilities, reached a contradiction
in each case, and so shown that $\Gamma$ must be non-singular for some
choice of $(w_1,w_2,w_3,w_4) \in Q'$.  This
then means that we couldn't have had so many quadruples in $Q'$,
and therefore $Q$.  Therefore, $|S+S| \gg n^2$, and we are done.

\section{Proof of Theorem \ref{main_theorem}}

Let us suppose that $\eps > 0$ is some constant that we
will allow to depend on $c$, $k$, and a certain parameter $t$
mentioned below, but is not allowed to depend on $\ell$.  This 
$\eps > 0$ will later be chosen small enough to make our proofs work.  Also,
we suppose that 
\begin{equation} \label{ellS}
|S^\ell|\ =\ |S.S...S|\ \leq\ n^c,
\end{equation}
where $\ell \geq 1$ is as large as we might happen to require, as a
function of $c$ and $k$ (and implicitly, $\eps$). 
\footnote{The reason we may choose $\ell$ as large as needed,
in terms of $c$ and $k$, is that we have freedom to choose 
$f(c,k,\ell)$ any way we please, so long as for fixed $c,k$ we have
$f(c,k,\ell) \to 0$ as $\ell \to \infty$.}

If (\ref{ellS}) holds, it follows that
$$
n = |S| \leq |S.S| \leq \cdots \leq |S^\ell| < n^c.
$$
Letting $M \geq 1$ be some integer depending on $c$ and $k$ that
we choose later, we have that for $\ell$ large enough,
it is obvious that for some $t < \ell/M$, 
\footnote{In other words, there must be a long interval $[t,Mt+1]$
such that for $j$ in this interval $S^j$ is not much smaller than $S^{Mt+1}$.} 
$$
|S^t|^{1+\eps}\ \geq\ |S^{Mt+1}|.
$$ 
Furthermore, if we let 
$$
S_M\ :=\ \{s^M\ :\ s \in S\},
$$
then since 
$$
S_M^t S = \{(s_1 \cdots s_t)^M s'\ :\ s_1,...,s_t,s' \in S\}\ 
\subseteq\ S^{Mt+1},
$$ 
it follows that
\begin{equation} \label{RS}
|S^t|\ =\ |S_M^t|\ \leq\ |S_M^t S|\ \leq\ |S^{Mt+1}|\ \leq\ |S^t|^{1+\eps}
\ =\ |S_M^t|^{1+\eps}.
\end{equation}
In other words, $S_M^t S$ is not appreciably larger than $S_M^t$. 

Let us now define
$$
R\ :=\ S_M^t\ =\ \{ (s_1 \cdots s_t)^M\ :\ s_i \in S\},
$$
so that (\ref{RS}) becomes
\begin{equation} \label{RS2}
|RS|\ \leq\ |R|^{1+\eps}.
\end{equation}

We now arrive at the following useful lemma.

\begin{lemma} \label{averaging_lemma} Suppose that 
(\ref{RS2}) holds.  Then, there exists 
$$
s \in S,\ {\rm and\ } r' \in R,
$$
such that for at least $n^{1-O(t \eps)}$ values $s' \in S$ we have that 
there exists $r \in R$ satisfying 
$$
rs\ =\ r' s'.
$$
\end{lemma}

\noindent {\bf Note.}  Since we can choose $\eps$ as small as desired
in terms of $t$, we have that $n^{1-O(t \eps)}$ is basically
$n^{1-\delta}$, where $\delta > 0$ is as small as we might happen to
require later on in the argument ($\delta$ is allowed to depend on
$c$ and $k$, but not on $\ell$).

\noindent {\bf Proof of the lemma.}  Inequality (\ref{RS2}) easily 
implies that
$$
|\{ (s,s',r,r') \in S^2 \times R^2\ :\ rs = r's'\}|\ \geq\ |R|^{1-\eps} n^2.
$$ 
So, extracting the pair $s,r'$ producing the maximal number of pairs
$s',r$ satisfying $rs = r's'$, we find that this pair leads to 
at least 
$$
n |R|^{-\eps}\ \geq\ n |S^t|^{-O(\eps)}\ \geq\ n^{1-O(t \eps)}
$$
such pairs $(s',r) \in S \times R$.   

\hfill $\blacksquare$
\bigskip

Let $S'$ denote the set of $s'$ produced by this lemma, for the fixed
pair $(r',s)$.  We will show that 
\begin{equation} \label{kS'}
|k S'|\ \geq\ n^{k - O(\delta)},
\end{equation}
from which it will follow that
$$
|k S|\ \geq\ |k S'|\ \geq\ n^{k-O(\delta)},
$$
thereby establishing Theorem \ref{main_theorem}, since the larger we
may take $\ell \geq 1$, the smaller we may take $\delta > 0$ 
(so, our $f(c,k,\ell) = O(\delta) \to 0$ for fixed $c,k$ as 
$\ell \to \infty$).

To prove (\ref{kS'}), it suffices to show that the only solutions to
\begin{equation} \label{nontrivial}
x_1 + \cdots + x_k\ =\ x_{k+1} + \cdots + x_{2k},\ x_i \in S',
\end{equation}
are trivial ones:  Suppose that, on the contrary, this equation has
a non-trivial solution.  Then, upon multiplying through by $r'$, we
are led to the equation
$$
sr_1 + \cdots + s r_k\ =\ s r_{k+1} + \cdots + sr_{2k},\ {\rm where\ } 
r_i \in R.
$$ 
Cancelling the $s$'s, and writing $r_i = y_i^M$, $y_i \in S^t$, produces
the equation  
\begin{equation} \label{yeqn}
y_1^M + \cdots + y_k^M\ =\ y_{k+1}^M + \cdots + y_{2k}^M. 
\end{equation}
In order to be able to apply Theorem \ref{poly_theorem} to this to
reach a contradiction, we need to show that this equation is non-trivial,
given that we have a non-trivial solution to (\ref{nontrivial}).
First observe that
$$
x_i = \lambda x_j\ \iff\ r' x_i = \lambda (r' x_j)\ \iff\ 
s r_i = \lambda s r_j\ \iff\ r_i = \lambda r_j,
$$  
which means that certain of the $y_i^M$'s are scalar multiples of
one another if and only if the corresponding $x_i$'s are scalar multiples
of one another, and in fact with the same scalars $\lambda$.

After cancelling common terms from both sides of (\ref{yeqn}), we 
move the remaining terms from among $y_{k+1}^M + \cdots + y_{2k}^M$
to the left-hand-side, writing
$$
-y_{k+j}^M\ =\ (e^{\pi i/M} y_{k+j})^M.
$$
We also must collect together duplicates into a
single polynomial:  Say, for example, $y_1 = \cdots = y_j$.  Then, we
collapse the sum $y_1^M + \cdots + y_j^M$ to $(j^{1/M} y_1)^M$.

For $M$ large enough in terms of $k$,  Theorem \ref{poly_theorem}
tells us that our collapsed equation can have no
non-trivial solutions, and therefore neither can (\ref{yeqn}).
Theorem \ref{main_theorem} is now proved.
 
\section{Proof of Theorem \ref{poly_theorem}}

The proof of this theorem will make use of the ideas that go into the
proof of the so-called $ABC$-theorem, which is also known as Mason's
Theorem \cite{mason} (see also \cite{granville} for a very nice
introduction).  Although there are versions of Mason's theorem already worked
out for some quite general contexts, in our useage of the ideas that go
into the proof of this theorem, we will need to allow some of the 
polynomials to have common factors.  Our proof is similar in many respects 
to the one appearing in \cite{bayat} and \cite{bayat2} --
the only difference is that we consider polynomials where none is a 
scalar multiple of another, and impose no coprimeness condition.

\subsection{The basic ABC theorem}

Before we embark on this task, let us recall the most basic ABC theorem,
and see its proof.

\begin{theorem} \label{ABC} Suppose that $A(x), B(x), C(x) \in \C[x]$ are coprime
polynomials, not all constant, such that 
$$
A(x) + B(x)\ =\ C(x).
$$
(Note that if any two share a common polynomial factor, then so must 
all three.)  Then, if we let $k$ denote the number of distinct roots of 
$A(x)B(x)C(x)$, we have that 
$$
{\rm max}({\rm deg}(A), {\rm deg}(B), {\rm deg}C)\ \leq\ k -1.
$$
\end{theorem} 

\noindent {\bf Remark 1.}  This theorem easily implies that 
the ``Fermat'' equation 
$$
f(x)^n + g(x)^n\ =\ h(x)^n
$$
has no solutions for $n \geq 3$, except trivial ones:  
Suppose that at most one of $f,g,h$ is 
constant, and that this equation does, in fact, have solutions.  Letting
$f$ be the polynomial of maximal degree, we find that
$$
n {\rm deg}(f)\ =\ {\rm deg}(f^n)\ \leq\ k - 1\ \leq\ {\rm deg}(fgh)-1
\ \leq\ 3 {\rm deg}(f)-1. 
$$
So, $n \leq 2$ and we are done.
\bigskip

\noindent {\bf Proof.}  The proof of the theorem makes use a remarkably
simple, yet powerful ``determinant trick''.  First, consider the 
determinant
\begin{equation} \label{the_matrix}
\Delta\ :=\ \left | \begin{array}{cc} A(x) & B(x) \\ A'(x) & B'(x) \end{array}
\right |.
\end{equation}
Note that this matrix is a Wronskian.

Let us see that $\Delta \neq 0$:  If $\Delta = 0$, we would have
that 
$$
A(x) B'(x)\ =\ B(x) A'(x).
$$
Since $A(x)$ and $B(x)$ are coprime, we must have that 
$$
A(x)\ |\ A'(x),\ {\rm and\ } B(x)\ |\ B'(x).
$$
Both of these are impossible, unless of course
both $A(x)$ and $B(x)$ are constants.  If both are constants, then
so is $C(x)$, and we contradict the hypotheses of the theorem.  So,
we are forced to have $\Delta \neq 0$.

Now suppose that 
$$
A(x) B(x) C(x)\ =\ c \prod_{i=1}^k (x - \alpha_i)^{a_i},
$$
so that $A,B,C$ have only the roots $\alpha_1,...,\alpha_k$, with
multiplicities $a_1,...,a_k$, respectively.  We will now see that
$$
R(x)\ :=\ \prod_{i=1}^k (x - \alpha_i)^{a_i-1}
$$
divides $\Delta$.  To see this, note that for each $i=1,...,k$,
$(x-\alpha_i)^{a_i}$ divides either $A(x),B(x),$ or $C(x)$, since
all three are coprime.  Note also that adding the first column  
of the matrix in (\ref{the_matrix})
to the second does not change the determinant, so that 
\begin{equation} \label{the_matrix2}
\Delta\ =\  \left | \begin{array}{cc} A(x) & C(x) \\ A'(x) & C'(x) 
\end{array} \right |.
\end{equation}
(Note that here we have used the fact that differentiation is a linear
map from the space of polynomials to itself.)  Now using the fact that
$$
(x-\alpha_i)^{a_i}\ |\ f(x)\ \ \Longrightarrow\ \ (x-\alpha_i)^{a_i-1}
\ |\ f'(x),
$$ 
it follows that $(x-\alpha_i)^{a_i-1}$ divides all the elements of 
some column of either the matrix (\ref{the_matrix}), or (\ref{the_matrix2}).
It follows that 
$$
(x-\alpha_i)^{a_i-1}\ |\ \Delta,\ {\rm and\ therefore\ }R(x) | \Delta.
$$
So,
$$
{\rm deg}(\Delta)\ \geq\ {\rm deg}(R(x))\ =\ {\rm deg}(ABC) - k.
$$
But we also have a simple upper bound on the degree of $\Delta$,
\begin{eqnarray*}
{\rm deg}(\Delta)\ &\leq&\ \min( {\rm deg}(AB) - 1, 
{\rm deg}(AC)-1, {\rm deg}(BC)-1) \\
&=&\ {\rm deg}(ABC) - 1 - \max({\rm deg}(A), {\rm deg}(B), {\rm deg}(C)).
\end{eqnarray*}
Collecting together the above inequalities clearly proves the theorem.
\hfill $\blacksquare$

\subsection{Two lemmas about polynomials}

We will require the following basic fact about Wronskians, which we
will not bother to prove.

\begin{lemma} \label{good_lemma}  
Suppose that $f_1,...,f_\ell \in \C[x]$, or even $\F_p[x]$, 
are non-zero polynomials and that
$$
\left | \begin{array}{cccc} f_1 & f_2 & \cdots & f_\ell \\ f_1' & f_2' & 
\cdots & f_\ell' \\ 
f_1'' & f_2'' & \cdots & f_\ell'' \\ \vdots & \vdots & \ddots & \vdots \\ 
f_1^{(\ell-1)} & f_2^{(\ell-1)} & \cdots & f_\ell^{(\ell-1)} 
\end{array} \right |\ =\ 0. 
$$
Then, we have that there are polynomials 
$$
\alpha_1, ..., \alpha_\ell,
$$
not all $0$, such that
\begin{equation} \label{c1}
\alpha_1 f_1 + \cdots + \alpha_\ell f_\ell\ =\ 0,
\end{equation}
where 
\begin{equation} \label{c2}
\alpha_1'\ =\ \cdots\ =\ \alpha_\ell'\ =\ 0.
\end{equation} 
Note that this last condition is equivalent to saying that
$\alpha_i$ are constants in the $\C[x]$ settting, and
$p$th powers of some polynomials in the $\F_p[x]$ setting. 
\end{lemma}

We will also require the following lemma.

\begin{lemma} \label{lemma2} 
For $k \geq 2$ and $\eps > 0$ there exists $M \geq 1$
such that the following holds:  Suppose that 
$$
f_1, ..., f_{k-1}\ \in\ \C[x]
$$
are linearly independent over $\C$.
Then, the equation 
$$
f_1^M + \cdots + f_{k-1}^M + G^M f_k\ =\ 0
$$
has no solutions with $G$ satisfying the following two conditions:
\bigskip

$\bullet$ ${\rm deg}(G)\ \geq\ \eps D$, where 
$$
D\ :=\ \max({\rm deg}(f_1), ..., {\rm deg}(f_{k-1}), {\rm deg}(G^M f_k)/M);
$$

$\bullet$ ${\rm deg}({\rm gcd}(G,f_j)) < \eps D/2k$, for all $j=1,...,k-1$.
\end{lemma}

\noindent {\bf Proof of the lemma.}  Suppose, in fact, that there are 
such solutions, where we will later choose $M$ purely as a function of 
$\eps$ and $k$, in order to reach a contradiction.

Then, consider the determinant
$$
\Delta\ :=\ \left | \begin{array}{ccc} f_1^M & \cdots & f_{k-1}^M \\
(f_1^M)' & \cdots & (f_{k-1}^M)' \\
\vdots & \ddots & \vdots \\
(f_1^M)^{(k-2)} & \cdots & (f_{k-1}^M)^{(k-2)} \end{array} \right |.
$$
The hypotheses of the lemma, along with Lemma \ref{good_lemma}, imply
that $\Delta \neq 0$.

It is easy to see, upon using your favorite expansion for the determinant,
that
$$
{\rm deg}(\Delta)\ \leq\ M {\rm deg}(f_1 \cdots f_{k-1}) - k + 2.
$$
On the other hand, the $j$th column of the matrix is divisble by 
$f_j^{M-k+2}$, so that
$$
(f_1 \cdots f_{k-1})^{M-k+2}\ |\ \Delta.
$$
And, upon adding the first $k-2$ columns to the last column, we
see that $\Delta$ is also divisible by $G^{M-k+2}$.  So,
$$
(f_1 \cdots f_{k-1})^{M-k+2} {G^{M-k+2} \over {\rm gcd}(G^{M-k+2},
f_1^M \cdots f_{k-1}^M)}\ |\ \Delta.
$$
Using again the hypotheses of the lemma, we can easily bound the degree
of this gcd from above by 
$$
M \sum_{j=1}^{k-2} \eps D/2k\ <\ \eps M D/2.
$$
So,
$$
{\rm deg}(\Delta)\ \geq\ (M-k+2) {\rm deg}(f_1 \cdots f_{k-1} G) - \eps M D/2.
$$
Combining our upper and lower bounds on ${\rm deg}(\Delta)$, we deduce
$$
{\rm deg}(G)\ \leq\ {k-2 \over M-k+2}({\rm deg}(f_1 \cdots f_{k-1})-1)
+ {\eps  M D \over 2 (M-k+2)}.
$$
This is impossible for $M$ large enough in terms of $\eps$ and $k$, since
$$
{\rm deg}(G)\ >\ \eps D,
$$
and our lemma is therefore proved.

\subsection{An iterative argument}

To apply the above lemmas, we begin by letting $\eps = 1$, and suppose
that 
\begin{equation} \label{solution}
f_1^M + \cdots + f_k^M\ =\ 0
\end{equation}
has a non-trivial solution (none is a scalar multiple of another), 
where $M$ can be taken as large as we might happen to need.  And,
assume that we have pulled out common factors among all the $f_i$.

We furthermore assume that any subset of size $k-1$ of the
polynomials $f_1^M,...,f_k^M$, is linearly independent over $\C$, since
otherwise we have a solution to (\ref{solution}), 
but with a smaller value of $k$ (we could pull the coefficients in the
linear combination into the $M$th powers, as every element of $\C$ has
an $M$th root in $\C$).

Without loss of generality, assume that $f_k$ has the highest degree
among $f_1,...,f_k$, and let $D$ denote its degree.

Now suppose that 
\begin{equation} \label{first_suppose}
{\rm deg}( {\rm gcd}(f_k, f_j))\ >\ \eps D/2k,\ {\rm for\ some\ } j=1,...,k-1.
\end{equation}
Without loss of generality, assume that $j=k-1$.
Then, we may write
$$
f_k^M + f_{k-1}^M\ =\ G_1^M g_k,\ {\rm where\ deg}(G_1) > \eps D/2k,
$$
so that
$$
f_1^M + \cdots + f_{k-2}^M + G_1^M g_k\ =\ 0.
$$
Note that this last term is non-zero by the assumption that no $k-1$
of the polynomials $f_i^M$ can be linearly dependent over $\C$.

Then, we set 
$$
\eps_1\ =\ \eps/2k,
$$
and observe that
$$
{\rm deg}(G_1)\ >\ \eps_1 D.
$$
On the other hand, if (\ref{first_suppose}) does not hold, then we
proceed on to the next subsection. 

Now suppose that
\begin{equation} \label{second_suppose}
{\rm deg}( {\rm gcd}(f_i, G_1))\ >\ \eps_1 D/2(k-1),\ {\rm for\ some\ }
i=1,...,k-2,
\end{equation}
where here we redefine $D$ to 
$$
D\ =\ \max({\rm deg}(f_1),...,{\rm deg}(f_{k-2}), {\rm deg}(G_1^M g_k)/M).
$$
(Note that we still have ${\rm deg}(G_1) > \eps_1 D$).
Without loss of generality, assume $i=k-2$.  Then, we may write
$$
f_k^M + f_{k-1}^M + f_{k-2}^M\ =\ G_2^M g_{k-2},\ 
{\rm where\ deg}(G_2) > \eps_1 D/2(k-1),
$$ 
so that
$$
f_1^M + \cdots + f_{k-3}^M + G_2^M g_{k-2}\ =\ 0.
$$
Of course, we also have to worry about wether this final
term is $0$, but that is not a problem as it would imply that 
$k-1$ of the polynomials $f_i^M$ are dependent over $\C$. 

On the other hand, if (\ref{second_suppose}) does not hold, then we
proceed on to the next subsection.

We repeat this process, producing $\eps_2, \eps_3, ...$, and 
$G_3, G_4,...$.  We cannot continue to the point where our equation is
$$
f_1^M + G_{k-2}^M g_2\ =\ 0,
$$
with $G_{k-2}$ non-constant, because then we would have that 
$f_1$ has a common factor with $G_{k-2}$,
and therefore all of $f_1,...,f_k$ would have to have a common factor.
So, the process must terminate before reaching $G_{k-2}$.

\subsection{Conclusion of the proof of Theorem \ref{poly_theorem}}

When we come out of the iterations in the previous section, we will
be left with an equation of the form
$$
f_1^M + \cdots + f_J^M + G_{k-J-1}^M g_{J-1}\ =\ 0,
$$
where 
$$
{\rm deg}(G_{k-J-1})\ >\ \gamma(k) D,
$$
where
$$
D\ =\ \max({\rm deg}(f_1), ..., {\rm deg}(f_J), 
{\rm deg}(G_{k-J-1}^M g_{J-1})/M)), 
$$ 
and where $\gamma(k)$ is a function depending only on $k$.  Applying now
Lemma \ref{lemma2} to this, we find that this is impossible once 
$M$ is large enough.  So, our theorem is proved. 

\section{Acknowledgment}

We would like to thank Jozsef Solymosi for his suggestion to 
include Theorem \ref{main_theorem3}.

\end{document}